\def\ver{Nov. 12, 2001, v.4}
\documentstyle{amsppt}
\magnification=1200
\hsize=6.5truein
\vsize=8.9truein
\topmatter
\title Cycle Map for Strictly Decomposable Cycles
\endtitle
\author Andreas Rosenschon and Morihiko Saito
\endauthor
\abstract We introduce a class of cycles, called nondegenerate, strictly
decomposable cycles, and show that the image of each cycle in this
class under the refined cycle map to an extension group in the derived
category of arithmetic mixed Hodge structures does not vanish.
This class contains certain cycles in the kernel of the Abel-Jacobi map.
The construction gives a refinement of Nori's
argument in the case of a self-product of a curve.
As an application, we show that a higher cycle which is not annihilated
by the reduced higher Abel-Jacobi map produces uncountably many
indecomposable higher cycles on the product with a variety having a
nonzero global 1-form.
\endabstract
\endtopmatter
\tolerance=1000
\baselineskip=12pt

\def\scirc{\raise.2ex\hbox{${\scriptstyle\circ}$}}
\def\ssb{\raise.2ex\hbox{${\scriptscriptstyle\bullet}$}}
\def\cD{{\Cal D}}
\def\cH{{\Cal H}}
\def\cO{{\Cal O}}
\def\cM{{\Cal M}}
\def\cP{{\Cal P}}
\def\bC{{\Bbb C}}

\def\bN{{\Bbb N}}
\def\bQ{{\Bbb Q}}
\def\bR{{\Bbb R}}
\def\bZ{{\Bbb Z}}
\def\txi{\widetilde\xi}

\def\oa{\overline{a}}
\def\ob{\overline{b}}
\def\oe{\overline{e}}
\def\ok{\overline{k}}
\def\oet{\overline{\eta}}
\def\ak{\langle k\rangle}

\def\div{\hbox{\rm div}}
\def\Hdg{\hbox{\rm Hdg}}
\def\Dec{\hbox{\rm Dec}\,}

\def\AJ{\text{\rm AJ}}
\def\an{\text{\rm an}}
\def\dec{\text{\rm dec}}
\def\ind{\text{\rm ind}}
\def\alg{\text{\rm alg}}
\def\hom{\text{\rm hom}}
\def\CH{\hbox{{\rm CH}}}
\def\Spec{\text{{\rm Spec}}\,}
\def\Ext{\hbox{{\rm Ext}}}
\def\Gr{\text{{\rm Gr}}}
\def\Im{\hbox{{\rm Im}}}
\def\Ker{\hbox{{\rm Ker}}}
\def\Hom{\hbox{{\rm Hom}}}
\def\codim{\hbox{{\rm codim}}}

\def\MHS{\text{{\rm MHS}}}
\def\HS{\text{{\rm HS}}}

\def\SameAuthor{\vrule height3pt depth-2.5pt width1cm}

\document
\centerline{\bf Introduction}

\bigskip\noindent
Let
$ X $ be a smooth projective complex variety, and let
$ \CH^{p}(X) $ be the Chow group of codimension
$ p $ cycles modulo rational equivalence.
Consider the subgroup
$ \CH_{\hom}^{p}(X) $ consisting of cycles homologically equivalent
to zero.
We have the Abel-Jacobi map to the intermediate Jacobian [24].
If
$ p > 1 $, it is known that this map is neither surjective nor
injective in general.
D. Mumford [30] showed that the kernel can be very large,
see also [7].
However it is not necessarily easy to construct a nonzero element
in the kernel, because the rational equivalence is rather
complicated for explicit calculations.
An argument constructing such an example was found by M. Nori in
the case of a self-product of a curve, see [38].
His argument involves the transcendental part of the second
cohomology, and requires some more calculations to verify the
hypothesis on cycles.
In this paper, we simplify and generalize it by introducing
the notion of a nondegenerate, strictly decomposable cycle.

Let
$ X_{1}, X_{2} $ be smooth projective varieties over
$ \bC $, and set
$ X = X_{1} \times X_{2} $.
Take divisors
$ \zeta_{i} \in \CH^{1}(X_{i}) $ which are homologically
equivalent to zero.
We assume that they are not torsion.
Consider
$ \zeta := \zeta_{1} \times \zeta_{2} $ in
$ \CH^{2}(X)_{\bQ} $.
This belongs to the kernel of the Abel-Jacobi map.
It may be zero, for example, if
$ X_{1} $ and
$ X_{2} $ are the same elliptic curve and
$ \zeta_{1} $ and
$ \zeta_{2} $ are the same cycle of the form
$ [P] - [Q] $, see (2.9.2).
So we are interested in the following question:
When is
$ \zeta $ not rationally equivalent to zero
(with rational coefficients)?
This is nontrivial even in the case of a self-product of an
elliptic curve.
In this paper we give a partial answer as follows.

\medskip\noindent
{\bf 0.1.~Theorem.} {\it
Assume there exists an algebraically closed
subfield
$ k $ of
$ \bC $ such that
$ X_{1} $,
$ X_{2} $ and
$ \zeta_{1} $ are defined over
$ k $, but
$ \zeta_{2} $ is not.
Then the cycle
$ \zeta $ does not vanish in
$ \CH^{2}(X)_{\bQ} $.
}

\medskip
This is an example of an exterior product of cycles [39] although
the varieties here are all defined over
$ k $.
In the curve case, Theorem (0.1) implies (by replacing
$ k $) that
$ ([P_{1}] - [Q_{1}]) \times ([P_{2}] - [Q_{2}]) $ does not vanish in
$ \CH^{2}(X)_{\bQ} $ if there exists an algebraically closed subfield
$ k $ of
$ \bC $ such that
$ X_{1} $,
$ X_{2} $,
$ Q_{1} $,
$ Q_{2} $ are defined over
$ k $, but
$ P_{1} $,
$ P_{2} $ are not, and if the field of definition
$ k(P_{1}) $ of
$ P_{1} $ is algebraically independent of
$ k(P_{2}) $ (i.e. if they generate a subfield of transcendence degree
$ 2 $).
Here
$ k(P_{i}) $ is the image of the morphism
$ k(X_{i}) \to \bC $ defined by
$ P_{i} $, see (2.9).

The proof of (0.1) uses a refined cycle map and the
associated Leray filtration on the Chow group, which follow
from the theory of (arithmetic) mixed sheaves [34], [36], [37]
(see also [1]) as conjectured in [4], [7].
This construction was originally considered in order to define an
analogue of ``spreading out" of algebraic cycles [7] for mixed Hodge
Modules (see [34], 1.9), and was greatly inspired by earlier work of
M.~Green [23], C.~Voisin [40], [41] and others.
We remark that our Leray filtration seems to coincide with a
filtration which has been studied recently by M.~Green in the case
the variety is defined over
$ k $ as explained
in his talk at Azumino, July 2000.
We can use also Green's formulation for the proof of Theorem (0.1),
see Remark (ii) of (1.5).
The theory of mixed Hodge modules is essential in the case the
variety is not defined over
$ k $.

Note that
$ \zeta_{1} $ in Theorem (0.1) can be a cycle of any codimension
provided that its image by the Abel-Jacobi map is not torsion.
Indeed, we have a more general assertion (2.6) whose proof is
reduced to a simple lemma on Hodge structures (2.7).
We can extend Theorem (2.6) to the case of higher Chow groups
$ \CH^{p}(X,m) $ [9] as in (3.3), which implies

\medskip\noindent
{\bf 0.2.~Theorem.} {\it
Let
$ \zeta = \zeta_{1}\times\zeta_{2}\in\CH^{p}(X,m)_{\bQ} $ with
$ \zeta_{1}\in\CH^{p_{1}}(X_{1})_{\bQ} $ and
$ \zeta_{2}\in\CH^{p_{2}}(X_{2},m)_{\bQ} $.
Then
$ \zeta $ does not vanish if there exists an algebraically closed
subfield
$ k $ of
$ \bC $ such that
$ X_{1} $,
$ X_{2} $ are defined over
$ k $ and the following condition is satisfied for
$ i = 1 $ or
$ 2 : $

\noindent
(i)
$ \zeta_{i} $ is defined over
$ k $, but
$ \zeta_{3-i} $ is not,

\noindent
(ii) the image of
$ \zeta_{i} $ by the cycle map (3.1.2) does not vanish,

\noindent
(iii)
$ p_{3-i} = 1 $.
}

\medskip
Here
$ \CH^{1}(X_{2},m) = 0 $ for
$ m > 1 $ and
$ \CH^{1}(X_{2},1) = \bC $ by loc.~cit.
So we may assume
$ X_{2} = pt $ if the above hypothesis is satisfied for
$ i = 1 $ and
$ m = 1 $.

Originally this theory was considered in the case when both
$ \xi_{1}^{1} $ and
$ \xi_{2}^{1} $ are nonzero with the notation of (2.1) in order to
show the nonvanishing of some extension classes.
In the study of Nori's construction [38], a special case of (2.6) was
obtained for a self-product of a curve of certain type [32].
Then these were generalized to (2.6).

As an application, we have the following.
For a smooth complex projective variety
$ X $, let
$ \CH_{\ind}^{p}(X,1)_{\bQ} $ denote the group of indecomposable higher
cycles with rational coefficients, see (4.1).
We have the reduced higher Abel-Jacobi map (4.1.2)
which is analogous to Griffiths' Abel-Jacobi map.

\medskip\noindent
{\bf 0.3.~Theorem.} {\it
Let
$ X = X_{1}\times X_{2} $ as above.
If
$ \Gamma(X_{1},\Omega_{X_{1}}^{1}) \ne 0 $ and the reduced higher
Abel-Jacobi map (4.1.2) for
$ X_{2} $ is not zero, then
$ \CH_{\ind}^{p+1}(X,1)_{\bQ} $ is uncountable.
More precisely, let
$ \zeta_{2} $ be an element of
$ \CH^{p}(X_{2},1) $ such that its image by (4.1.2) does
not vanish.
Let
$ \cP $ denote the Picard variety of
$ X_{1} $.
Assume
$ X_{i} $,
$ \zeta_{2} $ and
$ \cP $ are defined over
$ k $.
Then for any
$ \zeta_{1} \in \cP(\bC) \setminus \cP(\ok) $,
the image of
$ \zeta_{1}\times \zeta_{2} $ in
$ \CH_{\ind}^{p+1}(X,1)_{\bQ} $ is nonzero.
}

\medskip
This does not contradicts Voisin's conjecture [42]
on the countability of
$ \CH_{\ind}^{2}(X,1)_{\bQ} $, because
$ \CH_{\ind}^{1}(X_{2},1)_{\bQ} = 0 $.
It is known that
$ \CH_{\ind}^{p+1}(X,1)_{\bQ} $ is uncountable for certain varieties
$ X $, see [13].
The assumption on the nonvanishing of the image of
$ \zeta_{2} $ by (4.1.2) is satisfied for some examples, see [14], [35].

The paper is organized as follows.
In Sect.~1 we recall the notion of arithmetic mixed Hodge structures.
In Sect.~2 we prove Theorem (0.1).
Then Sect.~3 treats the case of higher cycles, and the application to
indecomposable higher cycles is given in Sect.~4.

In this paper, a variety means a separated scheme of finite type over
a field.
For a complex variety
$ X $ and a ring
$ A $, we denote
$ H^{j}(X^{\an},A) $ by
$ H^{j}(X,A) $ to simplify the notation.

\bigskip\bigskip
\centerline{{\bf 1.~Arithmetic mixed Hodge structure}}

\bigskip\noindent
{\bf 1.1.~Mixed Hodge structure.}
For a subfield
$ k $ of
$ \bC $,
let
$ \MHS_{k} $ denote the category of graded-polarizable mixed
$ \bQ $-Hodge structure with
$ k $-structure, see [17], [18], [26], etc.
(In this paper,
$ l $-adic cohomology is not used, because
$ k $ is assumed to be algebraically closed and the embedding of
$ k $ into
$ \bC $ is fixed.)
An object
$ H $ of
$ \MHS_{k} $ consists of a bifiltered
$ k $-vector space
$ (H_{k};F,W) $ and a filtered
$ \bQ $-vector space
$ (H_{\bQ},W) $ together with an filtered isomorphism
$ \alpha : (H_{k},W)\otimes_{k}\bC =
(H_{\bQ},W)\otimes_{\bQ}\bC $,
and they define a graded-polarizable mixed
$ \bQ $-Hodge structure in the sense of [16] after tensoring
$ H_{k} $ with
$ \bC $.
Here the polarizations on the graded pieces of
$ W $ are assumed to be defined over
$ k $ so that the graded pieces are semisimple.
Morphisms are pairs of morphisms of bifiltered or filtered vector
spaces compatible with
$ \alpha $.
This is an abelian category such that every morphism is bistrictly
compatible with
$ (F,W) $.
We have naturally a constant object
$ \bQ_{k} $.
We can define the Tate twist
$ (p) $ for
$ p \in \bZ $ as in loc. cit.

For
$ H, H' \in \MHS_{k} $,
we have
$$
\Ext_{\MHS_{k}}^{i}(H,H') = 0\quad \text{for }i > 1,
\leqno(1.1.1)
$$
by an argument similar to [12], see [36].

Let
$ X_{k} $ be a
$ k $-variety.
Then we have naturally the cohomology
$ H^{i}(X_{k},\bQ) $ in
$ \MHS_{k} $ by [16] together with the compatibility of de Rham
cohomology with base change.
Furthermore, we can define canonically
$ K_{\cH}(X_{k}) $ in the derived category
$ D^{b}\MHS_{k} $ such that its
cohomology
$ H^{i}K_{\cH}(X_{k}) $ is isomorphic to
$ H^{i}(X_{k},\bQ) $ in
$ \MHS_{k} $.
This can be done by applying the construction of the direct image
for perverse sheaves [5] to mixed Hodge Modules.
In the case
$ X_{k} $ is smooth and quasiprojective, we may assume that the
complement of
$ X_{k} $ in a smooth projective compactification of
$ X_{k} $ is a divisor.
Then we can take two sets of general hyperplane sections
$ \{H_{k,i}\} $ and
$ \{H'_{k,j}\} $ such that
$ \cap_{i} H_{k,i} = \cap_{j} H'_{k,j} = \emptyset $,
and use the Cech and co-Cech complexes associated to the affine
coverings
$ \{U_{k,i} = X_{k} \setminus H_{k,i}\} $ and
$ \{U'_{k,j} = X_{k} \setminus H'_{k,j}\} $ together with the
generalization of the weak Lefschetz theorem in [6], so that we get
a complex whose
$ p $-th component is
$$
\oplus_{|I|-|J|=p-\dim X}H^{\dim X}(U_{k,I},U_{k,I}\setminus
U'_{k,J};\bQ),
$$
where
$ U_{k,I} = \cap_{i\in I}U_{k,i} $,
$ U'_{k,J} = \cap_{j\in J}U'_{k,j} $.

By (1.1.1) we have a noncanonical isomorphism
$$
K_{\cH}(X_{k}) \simeq \oplus_{i} H^{i}(X_{k},\bQ)[-i]\quad
\text{in }D^{b}\MHS_{k}.
\leqno(1.1.2)
$$

\medskip\noindent
{\bf 1.2.~Deligne cohomology.}
We define an analogue of Deligne cohomology by
$$
H_{\cD}^{i}(X_{k},\bQ(j)) = \Ext_{{D}^{b}\MHS_{k}}^{i}
(\bQ_{k},K_{\cH}(X_{k})(j)).
$$
If
$ k = \bC $,
this coincides with the absolute
$ p $-Hodge cohomology of Beilinson [3] which is denoted by
$ H_{\cD}^{i}(X,\bQ(j)) $ in this paper, see also [20], [21], [25],
[33].
In general, if we put
$ X = X_{k}\otimes_{k}\bC $, then the forgetful functor induces
a natural morphism
$$
H_{\cD}^{i}(X_{k},\bQ(j)) \to H_{\cD}^{i}(X,\bQ(j)).
\leqno(1.2.1)
$$

More generally, for
$ H \in \MHS_{k} $,
we define
$$
H_{\cD}^{i}(X_{k},H(j)) =
\Ext_{{D}^{b}\MHS_{k}}^{i}(\bQ_{k},K_{\cH}(X_{k})\otimes H(j)).
$$
By (1.1.1) we have a natural short exact sequence
$$
\aligned
0 \to \Ext_{\MHS_{k}}^{1}(\bQ_{k},H^{i-1}(X_{k},H)(j))
&\to H_{\cD}^{i}(X_{k},H(j))
\\
&\quad\to\Hom_{\MHS_{k}}(\bQ_{k},H^{i}(X_{k},H)(j)) \to 0,
\endaligned
\leqno(1.2.2)
$$
where
$ H^{i}(X_{k},H) := H^{i}(X_{k},\bQ)\otimes H $.

Let
$ \CH^{p}(X_{k})_{\bQ} $ be the Chow group of codimension
$ p $ cycles with rational coefficients modulo rational equivalence.
Assume
$ X_{k} $ is smooth.
Then we have a cycle map (see e.g. [34])
$$
cl : \CH^{p}(X_{k})_{\bQ} \to H_{\cD}^{2p}(X_{k},\bQ(p)).
\leqno(1.2.3)
$$
This is compatible with the cycle map
$$
cl : \CH^{p}(X)_{\bQ} \to H_{\cD}^{2p}(X,\bQ(p))
\leqno(1.2.4)
$$
for
$ X = X_{k}\otimes_{k}\bC $ via (1.2.1).
The cycle map (1.2.4) induces Griffiths' Abel-Jacobi map [24]
tensored with
$ \bQ $:
$$
\CH_{\hom}^{p}(X)_{\bQ} \to J^{p}(X)_{\bQ} = \Ext_{\MHS}^{1}
(\bQ,H^{2p-1}(X,\bQ)(p)),
\leqno(1.2.5)
$$
where
$ \CH_{\hom}^{p}(X) $ denotes the subgroup consisting of
homologically equivalent to zero cycles, and the last isomorphism is
due to [12].

\medskip\noindent
{\bf 1.3.~Leray filtration.}
Let
$ S_{k} $ be a
$ k $-variety, and put
$ Y_{k} = X_{k}\times_{k}S_{k} $.
Then we have a canonical isomorphism
$$
K_{\cH}(Y_{k}) = K_{\cH}(X_{k})\otimes K_{\cH}(S_{k}).
\leqno(1.3.1)
$$
There is an increasing Leray filtration
$ F_{L} $ on
$ K_{\cH}(Y_{k}) $
induced by the canonical truncation
$ \{\tau_{\le r}\} $ on
$ K_{\cH}(X_{k}) $.
(In general,
$ \tau_{\le r}K^{i} $ for a complex
$ (K^{\ssb},d) $ is defined by
$ K^{i} $ for
$ i < r $,
$ \Ker\,d $ for
$ i = r $,
and
$ 0 $ for
$ i > r $ .)
We have a canonical isomorphism
$$
\Gr_{r}^{{F}_{L}}K_{\cH}(Y_{k}) = H^{r}(X_{k},\bQ)\otimes
K_{\cH}(S_{k})[-r].
\leqno(1.3.2)
$$

Since the filtration
$ F_{L} $ splits by (1.1.2), it induces a decreasing Leray
filtration
$ F_{L} $ on
$ H_{\cD}^{i}(Y_{k},\bQ(j)) $ such that
$$
\aligned
\Gr_{{F}_{L}}^{r}H_{\cD}^{i}(Y_{k},\bQ(j))
&= \Ext_{{D}^{b}\MHS_{k}}^{r}(\bQ_{k},H^{i-r}(X_{k},\bQ)\otimes
K_{\cH}(S_{k})(j))
\\
&= H_{\cD}^{r}(S_{k},H^{i-r}(X_{k},\bQ)(j)).
\endaligned
\leqno(1.3.3)
$$
Here the filtration
$ F_{L} $ is shifted as in [16].
We have a natural short exact sequence
$$
\aligned
0 \to \Ext_{\MHS_{k}}^{1}(\bQ_{k},H^{r-1}(S_{k},H^{i-r}
(X_{k},\bQ))(j))\to H_{\cD}^{r}(S_{k},H^{i-r}(X_{k},\bQ)(j))&
\\
\to \Hom_{\MHS_{k}}(\bQ_{k}, H^{r}(S_{k},H^{i-r}(X_{k},\bQ))(j))&
\to 0.
\endaligned
\leqno(1.3.4)
$$
The filtration
$ F_{L} $ and this short exact sequence are compatible with the
pull-back induced by the base change under a morphism of
$ k $-varieties
$ S'_{k} \to S_{k} $.

Assume
$ X $,
$ S $ are smooth.
By the cycle map (1.2.3) for
$ Y_{k} $,
we get the induced Leray filtration
$ F_{L} $ on
$ \CH^{p}(Y)_{\bQ} $ such that
$$
\Gr_{{F}_{L}}^{r}cl : \Gr_{{F}_{L}}^{r}\CH^{p}(Y_{k})_{\bQ}
\to \Gr_{{F}_{L}}^{r}H_{\cD}^{2p}(Y_{k},\bQ(p))
\leqno(1.3.5)
$$
is injective.

\medskip\noindent
{\bf 1.4.~Arithmetic mixed Hodge structure.}
Let
$ k $ be a subfield of
$ \bC $.
We assume
$ k $ is algebraically closed for simplicity (hence
$ l $-adic sheaves are not used in this paper).
Then the category of arithmetic mixed Hodge structures
$ \cM_{\ak} $ is defined to be the inductive limit of the
categories of mixed Hodge Modules on
$ S_{k} = \Spec R $ (i.e.
mixed Hodge Modules on
$ S := S_{k}\otimes_{k}\bC $ whose underlying bifiltered
$ \cD $-Modules and polarizations are defined over
$ S_{k}) $,
where
$ R $ runs over the finitely generated smooth
$ k $-subalgebra of
$ \bC $.
Although we can also define it to be the inductive limit of the
categories of admissible variations of mixed Hodge structures on
$ S_{k} $,
we need mixed Hodge Modules in order to calculate higher extension
groups, because the adjunction relation between direct images and
pull-backs does not hold for the derived category of admissible
variations.
There is a canonical pull-back functor
$$
\pi^{*} : \MHS_{k} \to \cM_{\ak}
\leqno(1.4.1)
$$
induced by the projection
$ S_{k} = \Spec R \to \Spec k $ for any finitely generated
smooth
$ k $-subalgebra
$ R $ of
$ \bC $.
Let
$ \bQ_{\ak} = \pi^{*}\bQ_{k} $.
We have also the forgetful functor
$$
\iota : \cM_{\ak} \to \MHS
\leqno(1.4.2)
$$
by restricting to the geometric generic point of
$ S_{k} $ defined by the inclusion
$ R \to \bC $.

For a complex algebraic variety
$ X $,
the cohomology
$ H^{i}(X,\bQ_{\ak}) $ is naturally defined in
$ \cM_{\ak} $ by taking an
$ R $-scheme
$ X_{R} $ such that
$ X_{R}\otimes_{R}\bC = X $,
because the direct image of the constant sheaf by
$ X_{R} \to \Spec R $ is defined by the theory of mixed Hodge
Modules.
Furthermore, it can be lifted to a complex
$ K_{\ak}(X) \in D^{b}\cM_{\ak} $ having a
canonical isomorphism
$$
H^{i}K_{\ak}(X) = H^{i}(X,\bQ_{\ak}).
$$
We define an analogue of Deligne cohomology by
$$
H_{\cD}^{i}(X,\bQ_{\ak}(j)) =
\Ext_{{D}^{b}\cM_{\ak}}^{i}(\bQ_{\ak},K_{\ak}(X)(j)).
$$

From now on, we assume
$ X = X_{k}\otimes_{k}\bC $ for a smooth
$ k $-variety
$ X_{k} $.
Then the argument is very much simplified, and we have natural
isomorphisms
$$
\aligned
&K_{\ak}(X) = \pi^{*}K_{\cH}(X_{k}),\quad
H^{i}(X,\bQ_{\ak}) = \pi^{*}H^{i}(X_{k},\bQ),
\\
&H_{\cD}^{i}(X,\bQ_{\ak}(j)) = \varinjlim\,
H_{\cD}^{i}(X_{k}\times_{k}S_{k},\bQ(j)).
\endaligned
\leqno(1.4.3)
$$
By (1.3) we have a Leray filtration
$ F_{L} $ on
$ K_{\ak}(X) $ such that
$$
\Gr_{{F}_{L}}^{r}H_{\cD}^{i}(X,\bQ_{\ak}(j)) = \varinjlim
\Ext_{{D}^{b}\MHS_{k}}^{r}(\bQ_{k},H^{i-r}
(X_{k},\bQ)\otimes K_{\cH}(S_{k})(j)).
\leqno(1.4.4)
$$

Since
$ \CH^{p}(X)_{\bQ} =
\varinjlim\,\CH^{p}(X_{k}\times_{k}S_{k})_{\bQ} $,
we have a cycle map
$$
cl : \CH^{p}(X)_{\bQ} \to H_{\cD}^{2p}(X,\bQ_{\ak}(p)),
\leqno(1.4.5)
$$
which induces the (decreasing) Leray filtration
$ F_{L} $ on
$ \CH^{p}(X)_{\bQ} $ such that
$$
\Gr_{{F}_{L}}^{r}cl : \Gr_{{F}_{L}}^{r}\CH^{p}(X)_{\bQ} \to
\Gr_{{F}_{L}}^{r}H_{\cD}^{2p}(X,\bQ_{\ak}(p))
\leqno(1.4.6)
$$
is injective.
This cycle map is compatible with (1.2.4).

\medskip\noindent
{\bf 1.5.~Remarks.} (i)
By the compatibility of the cycle maps (1.2.4) and (1.4.5),
$ {F}_{L}^{1}\CH^{p}(X)_{\bQ} $ coincides with the subgroup
$ \CH_{\hom}^{p}(X)_{\bQ} $ of homologically equivalent to zero
cycles, and
$ {F}_{L}^{2}\CH^{p}(X)_{\bQ} $ is contained in the kernel of
the Abel-Jacobi map
$ \CH_{\AJ}^{p}(X)_{\bQ} $.
Here we have equality if
$ p = \dim X $ (see e.g.
[36], 3.6).
This follows from Murre's result on Albanese motives [31] together
with the compatibility of the action of a correspondence.
Restricting to the subgroup
$ \CH_{\alg}^{p}(X)_{\bQ} $ of algebraically equivalent to
zero cycles, we have more generally the equality
$ {F}_{L}^{2}\CH_{\alg}^{p}(X)_{\bQ} = \CH_{\AJ}^{p}(X)
_{\bQ} \cap \CH_{\alg}^{p}(X)_{\bQ} $ for any
$ p $,
see [37], 3.9.

We can show
$ \Gr_{{F}_{L}}^{r}\CH^{p}(X)_{\bQ} = 0 $ for
$ r > p $.
In the case
$ k = \overline{\bQ} $,
it is expected that the filtration
$ F_{L} $ gives the conjectural filtration of Bloch [7] and
Beilinson [4], see also [27].
The problem is the injectivity of the cycle map (1.4.5).
In the case
$ p = 2 $ and
$ m = 0 $,
this can be reduced to a conjecture of Beilinson [4] and Bloch on
the injectivity of Abel-Jacobi map for codimension
$ 2 $ cycles defined over number fields (see [11] for a special case
where the conjecture is verified).

\medskip
(ii) It is also possible to consider a variant of Deligne
cohomology by replacing
$ X_{k}\times_{k}S_{k} $ with
$ X\times S $ in (1.4.3) (i.e. by forgetting
$ k $-structure).
Then we get a filtration similar to the one studied recently by
M. Green and also to the one treated in [28].
(However, the definition given in loc. cit. in terms of currents
does not seem to be correct even in the product case, because the
Hodge filtration
$ F $ on some graded pieces of
$ \Dec L $ on the complex of currents is filtered quasi-isomorphic
to the filtration
$ \sigma $ on the localization of holomorphic differential forms,
and not the pole order filtration of Griffiths [24] and Deligne [15].)
In view of the above conjecture on the injectivity of Abel-Jacobi
map, it may be expected that we get still the same filtration on
the Chow group after the modification.
Note that we can forget the
$ k $-structure for the proof of Theorems (0.1) and (0.2).

\medskip
(iii) For a subfield
$ K $ of
$ \bC $ containing
$ k $,
we can repeat the arguments in (1.4) by restricting
$ R $ to those contained in
$ K $.

\bigskip\bigskip
\centerline{{\bf 2.~Strictly decomposable cycles}}

\bigskip\noindent
{\bf 2.1.}
Let
$ X_{1} $ and
$ X_{2} $ be smooth complex projective varieties defined over an
algebraically closed subfield
$ k $ of
$ \bC $, i.e. there are smooth projective
$ k $-varieties
$ X_{i,k} $ such that
$ X_{i} = X_{i,k}\otimes_{k} \bC $.
Put
$ X_{k} = X_{1,k} \times_{k} X_{2,k} $,
$ X = X_{k} \otimes_{k} \bC $ as in the introduction.
We say that a cycle
$ \zeta $ on
$ X $ is strictly decomposable if there exist subfields
$ K_{i} $ of
$ \bC $ finitely generated over
$ k $,
together with cycles
$ \zeta_{i} $ on
$ X_{i,K_{i}} := X_{i} \otimes_{k} K_{i} $ for
$ i = 1, 2 $ such that the canonical morphism
$ K_{1} \otimes_{k} K_{2} \to \bC $ is injective, and
$ \zeta $ coincides with the base change of the cycle
$ \zeta_{1} \times_{k} \zeta_{2} $ on
$$
X_{1,K_{1}} \times_{k}X_{2,K_{2}} =
X_{k} \otimes_{k} (K_{1} \otimes_{k} K_{2})
$$
by
$ K_{1} \otimes_{k} K_{2} \to \bC $ (replacing
$ k $ with a finite extension if necessary).
We say that a strictly decomposable cycle is of bicodimension
$ (p_{1},p_{2}) $ if
$ \codim\, \zeta_{i} = p_{i} $.
Put
$ p = p_{1} + p_{2} $.

Let
$ R_{i} $ be a finitely generated smooth
$ k $-subalgebra of
$ K_{i} $ such that the fraction field of
$ R_{i} $ is
$ K_{i} $,
and
$ \zeta_{i} $ is defined over
$ R_{i} $.
Set
$ S_{i} = \Spec R_{i}\otimes_{k}\bC $.
We denote by
$ \zeta_{i}\otimes_{k}\bC $ the base change of
$ \zeta_{i} $ by
$ k \to \bC $.
Note that
$ \zeta_{i}\otimes_{k}\bC \in \CH^{p_{i}}(X_{i}\times S_{i}) $, and
its restriction to
$ X_{i}\times\{\oet_{i}\} $ coincides with
$ \zeta_{i}\otimes_{R_{i}}\bC $, where
$ \oet_{i} $ is defined by the inclusion
$ R_{i} \to \bC $.
Let
$$
\xi_{i}^{j} \in \Hom_{\MHS}(\bQ,H^{2p_{i}-j}(X_{i},\bQ)\otimes
H^{j}(S_{i},\bQ)(p_{i}))
\leqno(2.1.1)
$$
denote the K\"unneth components of the cycle class of
$ \zeta_{i}\otimes_{k}\bC $ in
$ H^{2p_{i}}(X_{i}\times S_{i},\bQ)(p_{i}) $.
Define
$$
H_{i} = H^{2p_{i}-1}(X_{i},\bQ)(p_{i}).
\leqno(2.1.2)
$$
If
$ \xi_{i}^{0} = 0 $,
then
$ \zeta_{i} $ belongs to
$ {F}_{L}^{1}\CH^{p_{i}}(X_{i,k} \otimes_{k} K_{i})_{\bQ} $
by (1.3) together with Remark (iii) in (1.5), and
$ \Gr_{{F}_{L}}^{1}cl(\zeta_{i}) $ gives
$$
\txi_{i}^{1} \in H_{\cD}^{1}(S_{i},H_{i})
\leqno(2.1.3)
$$
by forgetting the
$ k $-structure, see (1.3.3).

\medskip\noindent
{\bf 2.2.~Proposition.} {\it
With the above notation,
$ \xi_{i}^{1} = 0 $ if and only if
$ \txi_{i}^{1} $ comes from
$ \xi'_{i} \in \Ext_{\MHS}^{1}(\bQ, H_{i}) $ by the pull-back under
$ a_{S_{i}} : S_{i} \to \Spec \bC $.
In the case
$ p_{i} = 1 $, these conditions are further equivalent to that
$ \zeta_{i} $ comes from
$ \CH^{1}(X_{i,k})_{\bQ} $.
}

\medskip\noindent
{\it Proof.}
This follows from the short exact sequence (1.2.2) for
$ r = 1 $ and
$ k = \bC $, where the first term is isomorphic to
$ \Ext_{\MHS}^{1}(\bQ,H_{i}) $, and the first morphism is induced by
the pull-back functor for
$ a_{S_{i}} $.
Indeed,
$ \xi_{i}^{1} $ coincides with the image of
$ \txi_{i}^{1} $ in the last term, and we get the first equivalence.
Then the second is clear because the cycle map to Deligne cohomology is
an isomorphism in the divisor case and the base change by
$ k \to \bC $ induces an injective morphism of Chow groups with rational
coefficients.

\medskip\noindent
{\bf 2.3.~Remark.}
By (2.2) the cohomology class
$ \xi_{i}^{1} $ measures how the cycle
$ \zeta_{i}\otimes_{k}\bC $ varies along the parameter space
$ S_{i} $.
If the first equivalent conditions of (2.2) are satisfied, we may assume
$ R_{i} = k $ as long as
$ \Gr_{{F}_{L}}^{1}cl(\zeta_{i}) $ is concerned, e.g. if
$ p_{i} = 1 $.
(Indeed, we can restrict to a
$ k $-valued point of
$ \Spec R_{i} $.)

\medskip\noindent
{\bf 2.4.~Definition.} We say that
$ \zeta_{i} $ is {\it degenerate} if
$ \xi_{i}^{0} = \txi_{i}^{1} = 0 $.
A strictly decomposable cycle
$ \zeta $ is called {\it separately degenerate} if one of
$ \zeta_{i} $ is degenerate, and {\it cohomologically degenerate} if
$ \xi_{i}^{0} = \xi_{i}^{1} = 0 $ for both
$ i $,
and {\it degenerate} if it is separately degenerate or cohomologically
degenerate.
In other words,
$ \zeta $ is {\it nondegenerate} if one of
$ \xi_{i}^{0} $ and
$ \txi_{i}^{1} $ does not vanish for each
$ i $ and one of
$ \xi_{1}^{0} $,
$ \xi_{1}^{1} $,
$ \xi_{2}^{0} $,
$ \xi_{2}^{1} $ does not vanish.

\medskip\noindent
{\bf 2.5.~Remark.} Assume
$ p_{i} = 1 $.
Then
$ \zeta_{i} $ is degenerate if and only if it is zero in the Chow group
with rational coefficients.
If
$ \zeta $ is cohomologically degenerate, it comes from a cycle on
$ X_{k} $ by the pull-back under the projection to
$ X_{k} $, see (2.2).
If furthermore
$ k $ is a number field, it is expected that
$ \zeta $ is rationally equivalent to zero according to the conjecture
of Beilinson [2] and Bloch (see also [8], [26]).
If
$ \zeta $ is nondegenerate, however, we can detect it by the refined
cycle map due to the following:

\medskip\noindent
{\bf 2.6.~Theorem.} {\it
Let
$ \zeta $ be a nondegenerate, strictly decomposable cycle of codimension
$ p $, and
$ \zeta_{1} $,
$ \zeta_{2} $ be as in (2.1).
Then
$ cl(\zeta ) \ne 0 $.
More precisely, if
$ r = \#\{i : \xi_{i}^{0} = 0 \} $, then
$ \zeta \in {F}_{L}^{r}\CH^{p}(X)_{\bQ} $ and
$ \Gr_{{F}_{L}}^{r}cl(\zeta ) \ne 0 $.
}

\medskip\noindent
{\it Proof.}
Let
$ Y_{i,k} = X_{i,k}\times_{k} S_{i,k} $ and
$ Y_{k} = Y_{1,k}\times_{k} Y_{2,k} $.
Then the cycle map (1.2.3) gives
$$
\xi_{i} = cl(\zeta_{i}) \in H_{\cD}^{2p_{i}}(Y_{i},\bQ(p_{i})),\quad
\xi = cl(\zeta) \in H_{\cD}^{2p}(Y,\bQ(p)),
$$
and
$ \xi $ is the external product of
$ \xi_{1} $ and
$ \xi_{2} $ using (1.3.1).
Since
$ \xi_{i} \in F_{L}^{1}H_{\cD}^{2p_{i}}(Y_{i},\bQ(p_{i})) $ if
$ \xi_{i}^{0} = 0 $, we see that
$ \xi \in F_{L}^{r}H_{\cD}^{2p}(Y,\bQ(p)) $ and hence
$ \zeta \in {F}_{L}^{r}\CH^{p}(X)_{\bQ} $ for
$ r $ as above.
Then the assertion is clear if one of
$ \xi_{i}^{0} $ does not vanish.
Indeed, it is reduced to the case
$ R_{i} = k $ (replacing
$ k $ with a finite extension if necessary),
and follows from the nonvanishing of
$ \xi_{i'}^{0} $ or
$ \txi_{i'}^{1} $ for
$ i' = 3 - i $.
Thus we may assume
$ \xi_{1}^{0} = \xi_{2}^{0} = 0 $ and
$ r = 2 $.

By (1.3.1) we have the product
$$
\txi_{1}^{1}\otimes \txi_{2}^{1} \in H_{\cD}^{2}(S,H_{1}\otimes
H_{2})
\leqno(2.6.1)
$$
in the notation of (2.1.3), and it is induced by
$ \Gr_{F_{L}}^{2}\xi $.
We have to show that the pull-back of
$ \txi_{1}^{1}\otimes \txi_{2}^{1} $ by any dominant morphism
$ S' \to S $ does not vanish.
It is enough to show that the restriction of
$ \txi_{1}^{1}\otimes \txi_{2}^{1} $ to any nonempty open subvariety
$ U $ of
$ S $ does not vanish by the same argument as in [36], 2.6.
(Indeed, we can restrict to a subvariety of
$ S' $ which is finite \'etale over
an open subvariety of
$ S $,
and then take the direct image.)

If both
$ \xi_{1}^{1} $ and
$ \xi_{2}^{1} $ are nonzero, the assertion is easy, and follows from the
same argument as in [36], 4.4.
Indeed, if
$ H^{2}(S,\bQ)^{\vee} $ denotes the dual of
$ H^{2}(S,\bQ) $, then
$ \xi_{1}^{1}\otimes \xi_{2}^{1} $ corresponds to a morphism of mixed
Hodge structures
$$
\xi^{2} : H^{2}(S,\bQ)^{\vee} \to H^{2p-2}(X,\bQ)(p),
$$
and its image has level
$ 2 $ by hypothesis.
(Here the level of a Hodge structure is the difference between the maximal
and minimal numbers
$ p $ such that
$ \Gr_{F}^{p} \ne 0 $.)
This level does not change by replacing
$ S $ with any nonempty open subvariety
$ U $ of
$ S $, see [16].
So we get the assertion in this case.

Now it remains to consider the case
$ \xi_{1}^{1} \ne 0 $ and
$ \xi_{2}^{1} = 0 $ (hence
$ \xi'_{2} \ne 0) $.
By (2.2),
$ \txi_{1}^{1}\otimes \txi_{2}^{1} $ is the pull-back of
$ \txi_{1}^{1}\otimes \xi'_{2} $ by
$ S \to S_{1} $,
and we may replace
$ R_{2} $ with
$ k $.
So the assertion is reduced to the case
$ R_{2} = k $.
Then by (1.2.2) it is sufficient to show the nonvanishing of
$$
\xi_{1}^{1}\otimes \xi'_{2} \in
\Ext_{\MHS}^{1}(\bQ,H^{1}(S_{1},\bQ)\otimes H_{1}\otimes H_{2}),
$$
after replacing
$ S_{1} $ with any nonempty open subvariety of
$ S_{1} $.
Let
$$
H_{0} = H^{1}(S_{1},\bQ),\quad H'_{0} = \Gr_{1}^{W}H_{0},\quad
H''_{0} = \Gr_{2}^{W}H_{0}.
$$
Then
$ \xi_{1}^{1} \in \Hom_{\MHS}(\bQ, H_{0}\otimes H_{1}) $ comes from an
element of
$ u \in \Hom_{\MHS}(\bQ,H'_{0}\otimes H_{1}) $,
because
$ H_{1} $ is pure of weight
$ -1 $.
We denote by
$ e \in \Ext_{\MHS}^{1}(H''_{0},H'_{0}) $ the extension class associated
to the canonical short exact sequence.

Assume
$ \xi_{1}^{1}\otimes \xi'_{2} = 0 $.
Using the long exact sequence associated with
$$
0 \to H'_{0}\otimes H_{1}\otimes H_{2} \to
H_{0}\otimes H_{1}\otimes H_{2} \to
H''_{0}\otimes H_{1}\otimes H_{2} \to 0,
\leqno(2.6.2)
$$
this means that
$ u\otimes \xi'_{2} \in
\Ext_{\MHS}^{1}(\bQ,H'_{0}\otimes H_{1}\otimes H_{2}) $ coincides with
the composition of some
$ v \in \Hom_{\MHS}(\bQ,H''_{0}\otimes H_{1}\otimes H_{2}) $ with
$ e\otimes id $.
We have to show that this implies
$ u = v = 0 $.
Since
$ H''_{0} $ is isomorphic to a direct sum of
$ \bQ(-1) $,
$ e \in \Ext_{\MHS}^{1}(H''_{0},H'_{0}) $ is identified with a direct
sum of
$$
e_{j} \in \Ext_{\MHS}^{1}(\bQ(-1), H'_{0}),
$$
and
$ (e\otimes id)\scirc v $ with
$ \sum_{j} e_{j}\otimes v_{j} $,
where
$ v_{j} \in \Hom_{\MHS}(\bQ, H_{1}\otimes H_{2}(-1)) $.
In particular,
$ (e\otimes id)\scirc v $ factors through the inclusion
$ H'_{0}\otimes (\sum_{j} \Im\, v_{j}) \to
H'_{0}\otimes H_{1}\otimes H_{2} $.
On the other hand,
$ u\otimes \xi'_{2} $ factors through
$ (\Im\, u)\otimes H_{2} \to H'_{0}\otimes H_{1}\otimes H_{2} $.
Thus, by using a direct sum decomposition of
$ H'_{0}\otimes H_{1}\otimes H_{2} $ compatible with
$ (\Im\, u)\otimes H_{2} $ and
$ H'_{0}\otimes (\sum_{j} \Im\, v_{j}) $, the assertion is reduced to
the following:

\medskip\noindent
{\bf 2.7.~Lemma.} {\it
Let
$ H_{i} $ be a polarizable Hodge structure of weight
$ 1 $ for
$ i = 0, 1, 2 $.
Let
$ u \in \Hom_{\HS}(\bQ(-1),H_{0}\otimes H_{1}) $, and
$ v_{j} \in \Hom_{\HS}(\bQ(-1),H_{1}\otimes H_{2}) $.
Then}
$$
(\Im\, u)\otimes H_{2} \cap
H_{0}\otimes \bigl(\sum_{j} \Im\, v_{j}\bigr) = 0.
\leqno(2.7.1)
$$

\medskip\noindent
{\it Proof.}
It is sufficient to show the assertion for the underlying
$ \bR $-Hodge structure.
Using the projections to the direct factors, we may assume that
$ H_{0} $,
$ H_{1} $,
$ H_{2} $ are simple (i.e.
$ 2 $-dimensional), and furthermore
$ H_{0} $ is isomorphic to
$ H_{1}, H_{2} $ (because
$ u = 0 $ or
$ v_{j} = 0 $ otherwise).
Let
$ \{e, \oe\} $ be a basis of
$ H_{0} = H_{1} = H_{2} $ over
$ \bC $ such that
$ e $ generates
$ F^{1} $, and
$ \oe $ is the complex conjugate of
$ e $.
Then
$$
u = a\,e\otimes \oe + \oa \,\oe\otimes e\quad\text{and}\quad
v_{j} = b_{j}e\otimes \oe + \ob_{j}\oe\otimes e
$$
for
$ a, b_{j} \in \bC $.
In particular,
$ H_{0}\otimes (\sum_{j} \Im\, v_{j}) $ is contained in
the subspace generated by
$ H_{0}\otimes e\otimes\oe $ and
$ H_{0}\otimes\oe\otimes e $.
Thus we get (2.7.1).
This completes the proof of Theorem (2.6).

\medskip\noindent
{\bf 2.8.~Proof of Theorem (0.1).}
The assertion follows from (2.6) and (2.2).

\medskip\noindent
{\bf 2.9.~Examples.}
(i) Let
$ C_{1}, C_{2} $ be smooth proper curves over
$ \bC $,
and assume they are defined over
$ k $ (i.e.
$ C_{i} = C_{i,k}\otimes_{k}\bC $ for a curve
$ C_{i,k} $ over
$ k) $.
Let
$ P_{i}, Q_{i} \in C_{i}(\bC) $,
and assume that
$ Q_{1}, Q_{2} $ and
$ P_{2} $ are defined over
$ k $,
but
$ P_{1} $ is not, and that
$ [P_{2}] - [Q_{2}] $ is not torsion in the Jacobian.
Then Theorem (2.6) implies
$$
([P_{1}] - [Q_{1}]) \times ([P_{2}] - [Q_{2}]) \ne 0\quad
\text{in }\CH^{2}(C_{1}\times C_{2}).
\leqno(2.9.1)
$$
Note that the condition on
$ [P_{2}] - [Q_{2}] $ is satisfied if there exists an
algebraically closed subfield
$ k_{0} $ of
$ k $ such that
$ C_{2}, Q_{2} $ are defined over
$ k_{0} $,
but not
$ P_{2} $ (using the embedding of
$ C_{2} $ in the Jacobian).
Replacing
$ k $,
(2.9.1) holds if
$ Q_{1}, Q_{2} $ are defined over
$ k $,
and if the field of definition of
$ P_{1} $ (i.e. the image of
$ k(C_{1}) $ to
$ \bC $ by the morphism defined by
$ P_{1}) $ and that of
$ P_{2} $ are algebraically independent (i.e. if there exists an
algebraically closed subfield
$ k' $ of
$ \bC $ which contains
$ k $ and such that
$ P_{2} $ is defined over
$ k' $,
but
$ P_{1} $ is not).
See [1] for the case
$ C_{1} = C_{2}, P_{1} = P_{2} $ and
$ Q_{1} = Q_{2} $,
where it is assumed that the rank of the N\'eron-Severi group of
$ C_{1}\times C_{2} $ is
$ 3 $ and
$ \div(K_{C_{1}}) - (2g-2)[Q_{1}] $ is not torsion in the
Jacobian of
$ C_{1} $ (here
$ K_{C_{1}} $ is the canonical line bundle).
Note that the assertion does not hold for an elliptic curve by the
following.
\medskip
(ii) Let
$ E $ be an elliptic curve over
$ \bC $ with the origin
$ O $.
Then for any
$ P \in E(\bC) $
$$
([P] - [O]) \times ([P] - [O]) = 0
\quad\text{in}\,\, \CH^{2}(E\times E).
\leqno(2.9.2)
$$
This can be verified by taking
$ Q $ such that
$ 2Q = P $, and using the diagonal and antidiagonal curves passing
$ (Q,Q) $ (which is viewed as a new origin).

\bigskip\bigskip
\centerline{{\bf 3.~The higher cycle case}}

\bigskip\noindent
{\bf 3.1.}
For a smooth complex variety
$ X $ and a positive integer
$ m $,
let
$ \CH^{p}(X,m)_{\bQ} $ denote the higher Chow group of
$ X $ with rational coefficients [9].
By [36] we have the refined cycle map
$$
cl : \CH^{p}(X,m)_{\bQ} \to {H}_{\cD}^{2p-m}(X,\bQ_{\ak}(p)),
\leqno(3.1.1)
$$
which is compatible with the cycle map in [35]
$$
\CH^{p}(X,m)_{\bQ} \to {H}_{\cD}^{2p-m}(X,\bQ(p)).
\leqno(3.1.2)
$$
The latter is expected to be compatible with the one in [10], [19].
(For
$ m = 1 $, it is easy to verify this by reducing to the case
$ p = 1 $.)

We assume
$ X = X_{k}\otimes_{k}\bC $ and
$ X_{k} $ is smooth proper.
Then
$ \CH^{p}(X,m)_{\bQ} $ and
$ {H}_{\cD}^{2p-m}(X,\bQ_{\ak}(p)) $ have the Leray filtration
$ F_{L} $ as before, and (3.1.1) induces
$$
\aligned
\CH^{p}(X,m)_{\bQ}
&\to \Gr_{F_{L}}^{1}{H}_{\cD}^{2p-m}(X,\bQ_{\ak}(p))
\\
&=\Ext^{1}(\bQ_{\ak}, H^{2p-m-1}(X,\bQ_{\ak})(p)),
\endaligned
\leqno(3.1.3)
$$
because
$ H^{2p-m}(X,\bQ_{\ak}) $ is pure of weight
$ 2p - m $.
This is the inductive limit of
$$
\CH^{p}(X_{k}\times_{k}S_{k},m)_{\bQ} \to
{H}_{\cD}^{1}(S_{k},H^{2p-m-1}(X_{k},\bQ)(p)).
\leqno(3.1.4)
$$

Let
$ X_{i} $,
$ X_{i,k}\,\,(i = 1,2) $ and
$ X $,
$ X_{k} $ be as in (2.1).
In this paper we define the notion of a strictly decomposable higher cycle
(in the strong sense) by allowing a higher cycle for
$ \zeta_{2} $ in the notation of (2.1).
More precisely, we say that a higher cycle
$ \zeta \in \CH^{p}(X,m)_{\bQ} $ is {\it strictly decomposable}
if there exist subfields
$ K_{1}, K_{2} $ of
$ \bC $ and cycles
$ \zeta_{1} \in \CH^{p_{1}}(X_{1,k}\otimes_{k}K_{1})_{\bQ} $,
$ \zeta_{2} \in \CH^{p_{2}}(X_{2,k}\otimes_{k}K_{2},m)_{\bQ} $ such
that the
$ K_{i} $ satisfy the conditions in (2.1),
$ p_{1} + p_{2} = p $,
and
$ \zeta $ is the base change of
$ \zeta_{1}\times \zeta_{2} $ as before.

Let
$ \xi_{1}^{j} $,
$ \txi_{1}^{1} $ and
$ H_{1} $ be as in (2.1).
For
$ i = 2 $,
we have the K\"unneth component of the cohomology class
$$
\xi_{2}^{j} \in \Hom_{\MHS}(\bQ,H^{2p_{2}-m-j}(X_{2},\bQ)\otimes
H^{j}(S_{2},\bQ)(p_{2})).
$$
Note that
$ \xi_{2}^{0} = 0 $ in this case.
Let
$ H_{2} = H^{2p_{2}-m-1}(X_{2},\bQ)(p_{2}) $.
Then, applying (3.1.4) to
$ \zeta_{2} $ and forgetting the
$ k $-structure, we get
$$
\txi_{2}^{1} \in H_{\cD}^{1}(S_{2},H_{2}).
$$
By the exact sequence (1.2.2), it induces
$ \xi_{2}^{1} $, and comes from
$ \xi'_{2} \in \Ext_{\MHS}^{1}(\bQ, H_{2}) $ if
$ \xi_{2}^{1} = 0 $.
Thus we can define the notion of a nondegenerate, strictly decomposable
higher cycle in the same way as in (2.4).

\medskip\noindent
{\bf 3.2.~Remark.}
We have
$ m = 1 $ if
$ \xi_{2}^{1} \ne 0 $, see the proof of (3.3) below.

\medskip\noindent
{\bf 3.3.~Theorem.} {\it
The assertion of Theorem (2.6) holds also for nondegenerate, strictly
decomposable higher cycles
$ \zeta \in \CH^{p}(X,m)_{\bQ} $.
}

\medskip\noindent
{\it Proof.}
The argument is similar to (2.6).
We can identify
$ \xi_{2}^{1} $ with a morphism of mixed Hodge structures
$$
H^{1}(S_{2},\bQ)^{\vee} \to
H_{2}\,(= H^{2p_{2}-m-1}(X_{2},\bQ)(p_{2})),
\leqno(3.3.1)
$$
where the source is the dual of
$ H^{1}(S_{2},\bQ) $,
and the target is pure of weight
$ - m - 1 $.
Therefore the image of (3.3.1) is a direct sum of
$ \bQ(1) $ if
$ m = 1 $,
and vanishes otherwise.
For a nonempty open subvariety
$ S'_{2} $ of
$ S_{2} $,
we have the injectivity of the restriction morphism
$$
\Gr_{2}^{W}H^{1}(S_{2},\bQ) \to \Gr_{2}^{W}H^{1}(S'_{2},\bQ),
\leqno(3.3.2)
$$
because
$ H_{Z}^{1}(S_{2},\bQ) = 0 $ for
$ Z = S_{2}\setminus S'_{2} $.
Then, taking the tensor with
$ H^{2p_{1}}(X_{1},\bQ)(p_{1})\otimes H_{2} $, we get the
assertion in the case
$ \xi_{1}^{0} \ne 0 $ and
$ \xi_{2}^{1} \ne 0 $ (where we may assume
$ R_{1} = k $).
If
$ \xi_{1}^{0} \ne 0 $ and
$ \xi_{2}^{1} = 0 $, we may assume
$ R_{1} = R_{2} = k $, and the assertion is clear.
Thus the assertion is reduced to the case
$ \xi_{1}^{0} = 0 $.

The argument is similar if both
$ \xi_{1}^{1} $ and
$ \xi_{2}^{1} $ are nonzero (in particular,
$ m = 1) $.
For any nonempty open subset
$ S' $ of
$ S := S_{1}\times S_{2} $,
we have the injectivity of the restriction morphism
$$
\Gr_{3}^{W}H^{2}(S,\bQ) \to \Gr_{3}^{W}H^{2}(S',\bQ),
$$
because the local cohomology
$ {H}_{Z}^{2}(S,\bQ) $ is pure of type
$ (1,1) $, where
$ Z = S \setminus S' $.
Thus we get the assertion in this case, and we may assume
either
$ \xi_{1}^{1} $ or
$ \xi_{2}^{1} $ vanishes.

If
$ \xi_{1}^{1} \ne 0 $ and
$ \xi_{2}^{1} = 0 $ (hence
$ \xi'_{2} \ne 0) $,
then the assertion follows immediately from the long exact
sequence associated to (2.6.2), because
$ m > 0 $.
If
$ \xi_{1}^{1} = 0 $ and
$ \xi_{2}^{1} \ne 0 $ (hence
$ \xi'_{1} \ne 0 $ and
$ m = 1) $,
then we may assume
$ R_{1} = k $ and
$ H_{2} $ is a direct sum of
$ \bQ(1) $ because the image of (3.3.1) is a direct factor of the
target
(and the similar assertion holds over
$ k) $.
Furthermore, it is sufficient to show the nonvanishing of the composition
of
$ \xi'_{1}\otimes \xi_{2}^{1} $ with the morphism induced by the
projection
$ H^{1}(S_{2},\bQ) \to \Gr_{2}^{W}H^{1}(S_{2},\bQ) $ after
replacing
$ S_{2} $ with any nonempty open subvariety.
Thus the assertion follows from the injectivity of (3.3.2).
This finishes the proof of (3.3).

\medskip\noindent
{\bf 3.4.~Proof of Theorem (0.2).}
This follows from (3.3) and (2.2).

\medskip\noindent
{\bf 3.5.~Remarks.} (i)
Theorems (2.6) and (3.3) hold also for a variety
$ X $ over a subfield
$ K $ of
$ \bC $, where the
$ K_{i} $ are assumed to be subfields of
$ K $, see Remark (iii) of (1.5).

\medskip
(ii) Let
$ C $ be a curve of positive genus, and assume
$ X_{1} = X_{2} = C $,
$ K_{1} = k(X_{1,k}) $, and
$ K_{2} = k $.
Then (2.6) simplifies Nori's construction explained in [38]
(because (2.6) does not assume the condition on the relation to the
transcendental part of the second cohomology).
So Theorems (2.6) and (3.3) may be viewed as a refinement
of Nori's construction.

\newpage
\centerline{{\bf 4.~Application to indecomposable higher cycles}}

\bigskip\noindent
{\bf 4.1.}
Let
$ X $ be a smooth complex projective variety.
We have a canonical morphism
$$
\CH^{p-1}(X)_{\bQ} \otimes_{\bZ} \bC^{*} \to
\CH^{p}(X,1)_{\bQ},
\leqno(4.1.1)
$$
see [22], [29].
Its cokernel (resp. image) is denoted by
$ \CH_{\ind}^{p}(X,1)_{\bQ} $ (resp.
$ \CH_{\dec}^{p}(X,1)_{\bQ}) $.
An element of this group is called an indecomposable (resp. a
decomposable) higher cycle.
Let
$ \Hdg^{p-1} $ be the maximal subobject of
$ H^{2p-2}(X,\bQ) $ which is isomorphic to a direct sum of copies of
$ \bQ(1-p) $ (i.e. the subgroup of Hodge cycles).
It is a direct factor of
$ H^{2p-2}(X,\bQ) $ by semisimplicity.
The usual cycle map induces the reduced higher Abel-Jacobi map
$$
\CH_{\ind}^{p}(X,1)_{\bQ} \to
\Ext_{\MHS}^{1}(\bQ,(H^{2p-2}(X,\bQ)/\Hdg^{p-1})(p)),
\leqno(4.1.2)
$$
since the image of
$ \CH^{p-1}(X)_{\bQ} \otimes_{\bZ} \bC^{*} $ is contained
in
$ \Ext_{\MHS}^{1}(\bQ,\Hdg^{p-1}(p)), $ see loc.~cit.

\medskip\noindent
{\bf 4.2.~Theorem.} {\it
With the notation of (3.1) and (4.1), assume
$ m = 1 $,
$ \xi_{1}^{0} = 0 $,
$ \xi_{1}^{1} \ne 0, R_{2} = k $,
and the image of
$ \zeta_{2} $ by (4.1.2) for
$ X_{2} $ is nonzero.
Then
$ \zeta \ne 0 $ in
$ \CH_{\ind}^{p}(X,1)_{\bQ} $.
}

\medskip\noindent
{\it Proof.}
Assume
$ \zeta $ is decomposable, i.e.
$ \zeta $ coincides with the image of
$ \sum_{j} \zeta'_{j}\otimes \alpha_{j} $ where
$ \zeta'_{j} \in \CH^{p-1}(X)_{\bQ} $ and
$ \alpha_{j} \in \bC^{*} $.
Here we may assume that the
$ \alpha_{j} $ are linearly independent over
$ \bQ $ in
$ \bC^{*}\otimes_{\bZ}\bQ $,
and hence
$ \zeta'_{j} \in \CH_{\hom}^{p-1}(X)_{\bQ} $, because the image of
$ \zeta $ by the usual cycle map (3.1.2) vanishes due to the
assumption
$ \xi_{1}^{0} = 0 $
(indeed, the last assumption implies
$ r = 2 $ in (3.3)).
To simplify the notation, let
$ R = R_{1} $ and
$ S = \Spec R\otimes_{k}\bC $.
We may assume that the
$ \zeta'_{j} $ and
$ \alpha_{j} $ are defined over
$ R $ replacing
$ R $ if necessary.
We denote by
$ \zeta \otimes_{k}\bC $ the base change of
$ \zeta $ by
$ k \to \bC $ as before,
and similarly for
$ \zeta'_{j}\otimes_{k}\bC, $ etc.
We will induce a contradiction by showing that the image of
$ \zeta_{2} $ by (4.1.2) vanishes.
Since
$ \Hom_{\MHS}(\bQ, H^{2p_{2}-1}(X_{2},\bQ)(p_{2})) = 0 $,
we see that the image of
$ \zeta \otimes_{k}\bC $ in
$ \Hom_{\MHS}(\bQ, H^{2p-1}(S\times X,\bQ)(p)) $ vanishes, and so does
that of
$ \sum_{j} (\zeta'_{j}\otimes \alpha_{j})\otimes_{k}\bC $.

Let
$ H_{1} = H^{2p_{1}-1}(X_{1},\bQ)(p_{1}) $,
$ H_{2} = H^{2p_{2}-2}(X_{2},\bQ)(p_{2}) $ as before.
Let
$ \xi $ (resp.
$ \xi' $) denote the image of
$ \zeta \otimes_{k}\bC $ (resp.
$ \sum_{j} (\zeta'_{j}\otimes \alpha_{j})\otimes_{k}\bC $) in
$$
\Ext_{\MHS}^{1}(\bQ, H^{1}(S,\bQ)\otimes H_{1}\otimes H_{2}),
\leqno(4.2.1)
$$
using the cycle map (3.1.2) together with the K\"unneth decomposition.
We choose direct sum decompositions (see (1.1.2))
$$
K_{\cH}(X) \simeq \sum_{i} H^{i}(X,\bQ)[-i],\quad
K_{\cH}(S) \simeq \sum_{i} H^{i}(S,\bQ)[-i]\quad
\text{in }D^{b}\MHS.
\leqno(4.2.2)
$$
Applying this to the image of
$ \zeta'_{j}\otimes_{k}\bC $ by the cycle map to
$ {H}_{\cD}^{2p-2}(S\times X,\bQ(p-1)) $, and restricting to the
relevant K\"unneth component, we get
$$
\aligned
\xi'_{j,0}
&\in \Hom_{\MHS}(\bQ,H^{1}(S,\bQ)\otimes H_{1}\otimes H_{2}(-1)),
\\
\xi'_{j,1}
&\in \Ext_{\MHS}^{1}(\bQ,H_{1}\otimes H_{2}(-1)),
\endaligned
$$
because
$ H^{0}(S,\bQ) = \bQ $.
We see that the other K\"unneth components do not contribute to (4.2.1),
because the
$ \alpha_{j}\otimes_{k}\bC $,
which are viewed as sections of
$ {\cO}_{S}^{*} $,
induce
$$
\eta_{j,0} \in \Hom_{\MHS}(\bQ,H^{1}(S,\bQ)(1)),\quad
\eta_{j,1} \in \Ext_{\MHS}^{1}(\bQ,\bQ(1)),
$$
by using (4.2.2).
Thus
$$
\xi' =
\sum_{j} \eta_{j,0}\otimes \xi'_{j,1} +
\sum_{j} \eta_{j,1}\otimes \xi'_{j,0}.
$$
On the other hand,
$ \xi = \xi_{1}^{1}\otimes \zeta'_{2} $ with the notation of (3.2).

We have a direct sum decomposition
$$
H^{1}(S,\bQ) = H^{1}(S,\bQ)' \oplus H^{1}(S,\bQ)'',
\leqno(4.2.3)
$$
such that
$ H^{1}(S,\bQ)'' $ is a direct sum of copies of
$ \bQ(-1) $,
and
$$
\Hom_{\MHS}(\bQ(-1),H^{1}(S,\bQ)) =
\Hom_{\MHS}(\bQ(-1),H^{1}(S,\bQ)''),
$$
because the composition
$ \bQ(-1) \to H^{1}(S,\bQ) \to \Gr_{2}^{W}H^{1}(S,\bQ) $ splits by
semisimplicity.
Clearly
$ \eta_{j,0} \in \Hom_{\MHS}(\bQ,H^{1}(S,\bQ)''(1)) $,
and
$ \xi'_{j,0}, \xi_{1}^{1} $ come respectively from
$$
\aligned
\xi''_{j,0}
&\in \Hom_{\MHS}(\bQ,\Gr_{1}^{W}H^{1}(S,\bQ)'\otimes
H_{1}\otimes H_{2}(-1)),
\\
\xi''_{1}
&\in \Hom_{\MHS}(\bQ,\Gr_{1}^{W}H^{1}(S,\bQ)'\otimes H_{1}),
\endaligned
$$
because
$ \Gr_{1}^{W}H^{1}(S,\bQ)' = \Gr_{1}^{W}H^{1}(S,\bQ) $.
So we get
$ \sum_{j} \eta_{j,0}\otimes \xi'_{j,1} = 0 $ by the equality
$ \xi = \xi' $ together with the decomposition
(4.2.3).
Furthermore, we may replace
$ H^{1}(S,\bQ) $ with
$ \Gr_{1}^{W}H^{1}(S,\bQ) $,
using the long exact sequence associated to (2.6.2) as in the proof of
(3.3) (because
$ m > 0) $.
Thus the image of
$ \zeta_{2} $ by (4.1.2) vanishes by the following lemma which we
apply to
$$
\aligned
H_{1}
&= \Gr_{1}^{W}H^{1}(S,\bQ)\otimes H^{2p_{1}-1}(X_{1},\bQ)(p_{1}),
\\
H_{2}
&= (H^{2p_{2}-2}(X,\bQ)/\Hdg^{p_{2}-1})(p_{2}-1).
\endaligned
$$

\medskip\noindent
{\bf 4.3.~Lemma.} {\it
Let
$ H_{i} $ be a polarizable Hodge structure of weight
$ 0 $ for
$ i = 1, 2 $.
Assume
$ \Hom_{\MHS}(\bQ,H_{2}) = 0 $.
Let
$ u \in \Hom_{\MHS}(\bQ,H_{1}) $,
$ e \in \Ext_{\MHS}^{1}(\bQ,H_{2}(1)) $,
$ u_{j} \in \Hom_{\MHS}(\bQ,H_{1}\otimes H_{2}) $,
$ e_{j} \in \Ext_{\MHS}^{1}(\bQ,\bQ(1)) $ for
$ j = 1, \dots , r $.
If
$ u\otimes e = \sum_{j} u_{j}\otimes e_{j} $ in
$ \Ext_{\MHS}^{1}(\bQ,H_{1}\otimes H_{2}(1)) $, then
$ u\otimes e = 0 $.
}

\medskip\noindent
{\it Proof.}
Since
$ \Im\, u $ is a direct factor of
$ H_{1} $,
we may assume
$ H_{1} = \bQ $ by using the projection to the direct factor
$ H_{1} \to \bQ $.
Then
$ u_{j} = 0 $,
and the assertion is clear.
This completes the proof of (4.2).

\medskip\noindent
{\bf 4.4.~Proof of Theorem (0.3).}
The last assertion of (0.3) follows from (4.2) combined with (2.2).
If
$ \CH_{\ind}^{p+1}(X,1)_{\bQ} $ is countable, there are countably many
$ \zeta_{1,j} \in \cP(\bC) \,\,(j \in \bN) $ such that for any
$ \zeta_{1} \in \cP(\bC) $,
there exists
$ j $ such that the images of
$ \zeta_{1}\times \zeta_{2} $ and
$ \zeta_{1,j}\times \zeta_{2} $ in
$ \CH_{\ind}^{p+1}(X,1)_{\bQ} $ coincide.
We may assume
$ \zeta_{1,j} \in \cP(k) $ by replacing
$ k $ with a countably generated subfield of
$ \bC $.
Take any
$ \zeta_{1} \in \cP(\bC) \setminus \cP(\ok) $,
and let
$ \zeta_{1,j} \in \cP(k) $ such that the image of
$ (\zeta_{1} - \zeta_{1,j})\times \zeta_{2} $ in
$ \CH_{\ind}^{p+1}(X,1)_{\bQ} $ is zero.
But this contradicts the above assertion, because
$ \zeta_{1} - \zeta_{1,j} \notin \cP(\ok) $.
This terminates the proof of (0.3).

\medskip\noindent
{\bf 4.5.~Remark.}
It is known that
$ \CH_{\ind}^{p+1}(X,1)_{\bQ} $ is uncountable in some cases, see e.g.
[13] for the case of the Jacobian of a curve, and also [14], [35] for
examples such that the hypothesis of (0.3) is satisfied.
(We can show that a higher cycle on a product of elliptic curves in
[22] is decomposable by restricting to the generic fiber of the first
projection.)

\bigskip\bigskip\centerline{{\bf References}}\bigskip

\item{[1]}
Asakura, M., Motives and algebraic de Rham cohomology, in: The arithmetic
and geometry of algebraic cycles (Banff), CRM Proc. Lect. Notes, 24, AMS,
2000, pp. 133--154.

\item{[2]}
Beilinson, A., Higher regulators and values of
$ L $-functions, J. Soviet Math. 30 (1985), 2036--2070.

\item{[3]}
\SameAuthor, Notes on absolute Hodge cohomology, Contemporary Math. 55
(1986) 35--68.

\item{[4]}
\SameAuthor, Height pairing between algebraic cycles, Lect. Notes in
Math., vol. 1289, Springer, Berlin, 1987, pp. 1--26.

\item{[5]}
\SameAuthor, On the derived category of perverse sheaves, ibid.
pp. 27--41.

\item{[6]}
Beilinson, A., Bernstein, J. and Deligne, P., Faisceaux pervers,
Ast\'erisque, vol. 100, Soc. Math. France, Paris, 1982.

\item{[7]}
Bloch, S., Lectures on algebraic cycles, Duke University Mathematical
series 4, Durham, 1980.

\item{[8]}
\SameAuthor, Algebraic cycles and values of
$ L $-functions, J. Reine Angew. Math. 350 (1984), 94--108.

\item{[9]}
\SameAuthor, Algebraic cycles and higher
$ K $-theory, Advances in Math., 61 (1986), 267--304.

\item{[10]}
\SameAuthor, Algebraic cycles and the Beilinson conjectures,
Contemporary Math. 58 (1) (1986), 65--79.

\item{[11]}
Bloch, S., Kas, A. and Lieberman, D., Zero cycles on surfaces with
$ p_{g} = 0 $, Compos. Math. 33 (1976), 135--145.

\item{[12]}
Carlson, J., Extensions of mixed Hodge structures, in: Journ\'ees
de G\'eom\'etrie Alg\'ebrique d'Angers 1979, Sijthoff-Noordhoff
Alphen a/d Rijn, 1980, pp. 107--128.

\item{[13]}
Collino, A. and Fakhruddin, N., Indecomposable higher Chow cycles on
Jacobians, preprint.

\item{[14]}
del Angel, P. and M\"uller-Stach, S., The transcendental part of the
regulator map for
$ K_{1} $ on a mirror family of
$ K3 $ surfaces, preprint.

\item{[15]}
Deligne, P., Equations Diff\'erentielles \`a Points Singuliers
R\'eguliers, Lect. Notes in Math. vol.~163, Springer, Berlin, 1970.

\item{[16]}
\SameAuthor, Th\'eorie de Hodge I, Actes Congr\`es Intern. Math., 1970,
vol. 1, 425-430; II, Publ. Math. IHES, 40 (1971), 5--57; III ibid., 44
(1974), 5--77.

\item{[17]}
\SameAuthor, Valeurs de fonctions
$ L $ et p\'eriodes d'int\'egrales,
Proc. Symp. in pure Math., 33 (1979) part 2, pp. 313--346.

\item{[18]}
Deligne, P., Milne, J., Ogus, A. and Shih, K., Hodge Cycles, Motives, and
Shimura varieties, Lect. Notes in Math., vol 900, Springer, Berlin, 1982.

\item{[19]}
Deninger, C. and Scholl, A., The Beilinson conjectures, Proceedings
Cambridge Math. Soc. (eds. Coats and Taylor) 153 (1992), 173--209.

\item{[20]}
El Zein, F. and Zucker, S., Extendability of normal functions associated
to algebraic cycles, in: Topics in transcendental algebraic geometry, Ann.
Math. Stud., 106, Princeton Univ. Press, Princeton, N.J., 1984, pp.
269--288.

\item{[21]}
Esnault, H. and Viehweg, E., Deligne-Beilinson cohomology, in:
Beilinson's conjectures on Special Values of
$ L $-functions, Academic Press, Boston,
1988, pp. 43--92.

\item{[22]}
Gordon, B. and Lewis, J., Indecomposable higher Chow cycles on products
of elliptic curves, J. Alg. Geom. 8 (1999), 543--567.

\item{[23]}
Green, M., Griffiths' infinitesimal invariant and the Abel-Jacobi
map, J. Diff. Geom. 29 (1989), 545--555.

\item{[24]}
Griffiths, P., On the period of certain rational integrals I, II, Ann.
Math. 90 (1969), 460--541.

\item{[25]}
Jannsen, U., Deligne homology, Hodge
$ D $-conjecture, and motives, in Beilinson's conjectures on Special
Values of
$ L $-functions, Academic Press, Boston, 1988, pp. 305--372.

\item{[26]}
\SameAuthor, Mixed motives and algebraic
$ K $-theory, Lect. Notes in Math., vol. 1400, Springer, Berlin, 1990.

\item{[27]}
\SameAuthor, Motivic sheaves and filtrations on Chow groups, Proc. Symp.
Pure Math. 55 (1994), Part 1, pp. 245--302.

\item{[28]}
Lewis, J.D., A filtration on the Chow groups of a complex projective
variety, Compos. Math. 128 (2001), 299--322.

\item{[29]}
M\"uller-Stach, S., Constructing indecomposable motivic cohomology
classes on algebraic surfaces, J. Alg. Geom. 6 (1997), 513--543.

\item{[30]}
Mumford, D., Rational equivalence of
$ 0 $-cycles on surfaces, J. Math. Kyoto Univ. 9 (1969), 195--204.

\item{[31]}
Murre, J.P., On the motive of an algebraic surface, J. Reine Angew. Math.
409 (1990), 190--204.

\item{[32]}
Rosenschon, A., A Note on arithmetic mixed Hodge structures,
preprint 2001.

\item{[33]}
Saito, M., Mixed Hodge Modules, Publ. RIMS, Kyoto Univ., 26 (1990),
221--333.

\item{[34]}
\SameAuthor, On the formalism of mixed sheaves, RIMS-preprint 784, Aug.
1991.

\item{[35]}
\SameAuthor, Bloch's conjecture, Deligne cohomology and higher Chow
groups, preprint RIMS-1284 (or math.AG/9910113).

\item{[36]}
\SameAuthor, Arithmetic mixed sheaves, Inv. Math. 144 (2001), 533--569.

\item{[37]}
\SameAuthor, Refined cycle maps, preprint (math.AG/0103116).

\item{[38]}
Schoen, C., Zero cycles modulo rational equivalence for some varieties
over fields of transcendence degree one, Proc. Symp. Pure Math. 46 (1987),
part 2, pp. 463--473.

\item{[39]}
\SameAuthor, On certain exterior product maps of Chow groups, Math. Res.
Let. 7 (2000), 177--194,

\item{[40]}
Voisin, C., Variations de structures de Hodge et z\'ero-cycles sur
les surfaces g\'en\'erales, Math. Ann. 299 (1994), 77--103.

\item{[41]}
\SameAuthor, Transcendental methods in the study of algebraic cycles,
Lect. Notes in Math. vol. 1594, pp. 153--222.

\item{[42]}
\SameAuthor, Remarks on zero-cycles of self-products of varieties, in:
Moduli of Vector Bundles, Lect. Notes in Pure and Applied Mathematics,
vol. 179, M. Dekker, New York, 1996, pp. 265--285.

\bigskip\noindent
Andreas Rosenschon

\noindent
Department of Mathematics

\noindent
Duke University

\noindent
Durham, NC 27708

\noindent
U.S.A

\noindent
E-Mail: axr\@math.duke.edu

\medskip\noindent
Morihiko Saito

\noindent
RIMS Kyoto University,

\noindent
Kyoto 606--8502

\noindent
Japan

\noindent
E-Mail: msaito\@kurims.kyoto-u.ac.jp

\medskip\noindent
\ver

\bye